\documentclass[11pt]{amsart}

\newtheorem{theorem}{Theorem}[section]
\newtheorem{lemma}[theorem]{Lemma}
\newtheorem{corollary}[theorem]{Corollary}
\newtheorem{proposition}[theorem]{Proposition}

\theoremstyle{definition}

\theoremstyle{remark}

\newcommand{\R}{{\mathbb R}}
\newcommand{\C}{{\mathbb C}}
\newcommand{\cD}{{\mathcal D}}

\begin{document}

\title[Asymptotics for Partial Integrodifferential Equations]
{Asymptotic Self-Similarity for Solutions of Partial Integrodifferential
Equations}

\author{Hans Engler}
\address{Dept. of Mathematics, Georgetown University\\
Box 571233\\ Washington, DC 20057\\
USA}
\email{engler@georgetown.edu} \maketitle

\begin{center}
{\bf Abstract}
\end{center}

\noindent

{\footnotesize The question is studied whether weak solutions of
linear partial integrodifferential equations approach a constant
spatial profile after rescaling, as time goes to infinity. The
possible limits and corresponding scaling functions are identified
and are shown to actually occur. The limiting equations are
fractional diffusion equations which are known to have
self-similar fundamental solutions. For an important special case,
is is shown that the asymptotic profile is Gaussian and
convergence holds in $L^2$, that is, solutions behave like
fundamental solutions of the heat equation to leading order.
Systems of integrodifferential equations occurring in
viscoelasticity are also discussed, and their solutions are shown
to behave like fundamental solutions of a related Stokes system.
The main assumption is that the integral kernel in the equation is
regularly varying in the sense of Karamata.}


\date{October 9, 2005}

\section{Introduction}
Consider the linear heat equation $u_t = \Delta u$ in $\R^n$, with
fundamental solution $U(x,t) = \frac{1}{\left(4\pi t\right)^{n/2}}
e^{-|x|^2/4t}$. It is well-known and easy to see from the solution
formula that as $t \to \infty$, $u(x,t) = U_0 U(x,t) +
o(t^{-n/2})$, where $U_0 = \int_{\R^n} u(\cdot,0)$ is the initial
mass of the solution, assumed to be finite. Thus, $t^{n/2}
u(x\sqrt{t},t) \to U_0 U(x,1)$. Similarly, for solutions of the
wave equation $u_{tt} = u_{xx}$ on $\R \times [0,\infty)$ with
initial data $u(\cdot,0) = u_0, \, u_t(\cdot,0) = 0$, the
well-known solution formula $u(x,t) = \frac{1}{2}\left(u_0(x+t) +
u_0(x-t)\right)$ implies that $tu(xt,t) \sim U_0
\frac{1}{2}\left(\delta_{-1} + \delta_1\right)$ as $t \to \infty$,
in this case in the sense of distributions. For the case of the
wave equation on $\R^n$, the solution formulae that use spherical
means imply that $t^n u(xt,t) \to U_0w$, where $w$ is a
distribution of dimension dependent order that is supported on the
unit sphere in $\R^n$. In the case of the heat equation, the
solution depends (up to a multiplicative factor) asymptotically on
the similarity variable $\xi = \frac{|x|}{\sqrt{t}}$, and
convergence is uniform. In the case of the wave equation, the
similarity variable is $\xi = \frac{|x|}{t}$, with convergence in
a space of distributions.

\medskip
In this paper, I investigate whether solutions of
integrodifferential equations
\[
u_t(\cdot,t) = a_0 \Delta u(\cdot,t) + \int_0^t a(t-s) \Delta u(
\cdot,s) ds
\]
in $\R^n$ have similar properties. Here $\R \ni a_0 \ge 0$ is a
scalar, and $a:[0,\infty) \to \R$ is a scalar kernel. The special
cases $a_0 = 1, \, a = 0$ and $a_0 = 0, \, a = 1$ correspond to
the heat equation and the wave equation, respectively, and
therefore are included. Thus the question is whether
$m(t)u(k(t)x,t) \sim u_\infty$ as $ t \to \infty$, in a suitable
sense.

\medskip
Since all equations in this class are of the form $u_t + \nabla
\cdot q = 0$ for some flux $q$, the $L^1$-integral of solutions is
formally preserved as $t$ varies,
\[\int_{\R^n} u(\cdot,t) = \int_{\R^n} u_0 = U_0 \, .
\]
One therefore expects that $m(t) = k(t)^n$, that is
\[
k(t)^n u(x k(t),t) \sim U_0w_\infty \] for a suitable function
$k(\cdot)$ as $t \to \infty$ in a suitable distributional sense,
for a limiting distribution $w_\infty$. Then the "trivial"
behavior $w_\infty = \delta_0$ can always be achieved by letting
$k$ grow to $\infty$ very rapidly. This trivial behavior must
therefore be excluded. Also, solutions are expected to go to zero
locally, so the trivial case $w_\infty = 0$ is possible if $k$
goes to $\infty$ too slowly and must also be excluded. With the
right choice of $k$ one hopes to obtain a nontrivial limit
$w\infty$.

\medskip
In this paper, I show that for a large class of such
integrodifferential equations there is a choice of $k$ (unique up
to an asymptotically constant factor) for which the limit
$w_\infty$ is indeed non-trivial. It turns out that the correct
choice is
\[
k(t) = \sqrt{t \left(a_0 + \int_0^t a(s)ds\right)} \, .
\]
The limiting distributions $w_\infty$ are also identified. They
turn out to belong to a one parameter family, parametrized by
$\beta \in (-1, 1]$, with $\beta = 0$ corresponding to the heat
equation, $\beta = 1$ corresponding to the wave equation, and the
cases of non-integer $\beta$ corresponding to fundamental
solutions of {\it fractional diffusion} equations. The main
assumption is that the integrated kernel $A(t) = a_0 + \int_0^t
a(s)ds$ should be regularly varying in the sense of
\textsc{Karamata} (\cite{bingham}), and the index of variation
$\beta$ then determines the limiting distribution $w_\infty$. As
an aside, it should be noted that for the same $w_\infty$, there
are many types of scaling functions $k$ possible that are not
asymptotically equivalent. It will be shown that all these
possible limiting distributions are actually attained (in the
sense of distribution, or in an important special case in $L^2$).

\medskip
The literature on self-similar asymptotics is huge, so I only
mention the book \cite{barenblatt} by G. Barenblatt. Fractional
diffusion equations were discussed  in \cite{schneiderwyss} and
\cite{fujita}, with systematic studies carried out in
\cite{eidelman, mainardi1, mainardi0, mainardi2}. Various physical
models leading to fractional diffusion equations are discussed in
\cite{hanyga1, hanyga2, zaslavsky}. The main reference for
integrodifferential equations of the type discussed here is the
book \cite{pruss} by J. Pr\"uss. The idea that regularly varying
integral kernels lead to asymptotically self-similar wave profiles
for problems in viscoelasticity is exploited in \cite{pipkin}, for
the case of the signalling problem. A related asymptotic concept
is equipartition of energy, discussed for a class of exponential
kernels corresponding to $\beta = 1$ in \cite{dassios}.

\medskip
The plan of this paper is as follows. In the following section 2,
two types of scaling (introduced at the end of this section) and
the relations between them are discussed. All possible
distributional limits and their corresponding scaling functions
$k$ are identified in section 3. Section 4 is devoted to giving
sufficient conditions such that these distributional limits are
actually attained. In section 5, the question of asymptotic
self-similarity is studied in $L^2$, leading to results about the
time-asymptotic behavior of solutions that are sharp to leading
order. The same question is taken up for three-dimensional linear
homogeneous isotropic viscoelasticity in section 6. Appendix A
contains two important technical results for families for scalar
integral equations, and two explicit examples for asymptotic
behavior outside the theory developed in this paper are presented
in appendix B.

\medskip
The notation $\langle u, \varphi \rangle$ will be used for the
result of applying a distribution $u \in \cD'$ to a test function
$\varphi \in \cD = C_0^\infty(\R^n)$. The pairing between test
functions $\Phi \in C_0^\infty(\R^n \times [0,\infty))$ and
distributions $U$ on $\R^n \times [0,\infty)$ is denoted by
$\langle \langle U,\Phi\rangle \rangle$. In particular, if
$U(\cdot,t)_{t \ge 0}$ is a family of distributions in $\cD'$ that
is measurable and bounded with respect to the system of seminorms
defining the usual topology on $\cD'$, one can write
\[
\langle \langle U,\Phi\rangle \rangle = \int_0^\infty \langle
U(\cdot,t), \Phi(\cdot,t) \rangle dt \, .
\]
Convolution with respect to $t \in \R$ is denoted by an asterisk,
$u*v(t) = \int_0^tu(t-s)v(s) \, ds$ if $u$ and $v$ are both
supported on the positive half axis. The Fourier transform of a
function $f \in \cD $ is denoted by $\hat f, \, \hat f(\xi) =
\int_{\R^n} e^{-i x \xi} f(x) \, dx$, and this is extended in the
usual way to functions in $L^1$ or in $L^2$ or to distributions.
The Laplace transform of a function $a:[0,\infty) \to \R$ is
denoted by $\tilde a(s) = \int_0^\infty e^{-st} a(t) dt$ if
defined, i.e. for $s \in \C$ such that $\Re s > \alpha$ for some
$\alpha \in (-\infty, +\infty]$. The usual notation for Lebesgue
spaces $L^p(\Omega)$ and for the Sobolev spaces $H^s(\R^n)$ is
employed, $1 \le p \le \infty, \, -\infty < s < \infty$.
Vector-valued Lebesgue and Sobolev spaces are denoted in the usual
way, e.g. $L^p([0,T],H^1(\R^n))$ or $H^s(\R^3,\R^3)$.

Let $\R^+ \ni t \mapsto u(\cdot,t)$ be a measurable and locally
bounded family of distributions on $\R^n$, and let $k:[0,\infty)
\to \R^+$ be continuous and increasing to $\infty$. In this paper,
the scaled version of $u$ is denoted by $u_k$, defined by
\[\langle u_k(\cdot,t),
\varphi \rangle = \langle u(\cdot,t), \varphi
(k(t)^{-1}\cdot)\rangle \] or, in case $u(\cdot,t)$ is a function
for almost all $t$,
\[  u_k(x,t) = k(t)^nu(k(t)x,t) \, .\]
Thus if $u(\cdot,t) \in L^1(\R^n)$, then also $u_k(\cdot,t) \in
L^1(\R^n)$, and the $L^1$ integral is unchanged.

There is an alternative scaling $u_{T,k}$, defined by
$$\langle u_{T,k}(\cdot,t), \varphi\rangle = \langle u(\cdot, Tt),
\varphi(k(T)^{-1} \cdot)\rangle$$ or if $u(\cdot, t)$ is a
function for almost all $t$,
\[u_{T,k}(x,t) = k(T)^nu(k(T)x,Tt) \]
for $x \in \R^n$ and $ t > 0$. The result now depends on $k$ and
the parameter $T > 0$. This scaling again preserves the
$L^1$-integral.

\section{Scaling}
Let $t \mapsto u(\cdot,t)$ be a measurable and locally bounded
(with respect to the usual system of seminorms) family of tempered
distributions on $\R^n$, and let $k:[0,\infty) \to \R^+$ be
continuous and increasing to $\infty$. Assume that the scaled
family $u_k(\cdot,t)$ converges in $\cD'$ to some $U \in \cD'$ as
$t \to \infty$, that is,
\begin{equation}
\label{one} \langle u_k(\cdot,t), \varphi(\cdot)\rangle \to \,
\langle U, \varphi \rangle \, \forall \varphi \in C^\infty_0(\R^n)
\, .
\end{equation}
It is necessary to relate this property to the behavior of the
family $u_{T,k}$ as $T \to \infty$. Using a tensor product
argument, one easily sees that for all test functions $\Phi \in
C^\infty_0(\R^n \times (0,\infty))$
\begin{equation}
\label{two} \int_0^\infty \langle  u_{T,k}(\cdot,\tau),
\Phi(\cdot,\tau)\rangle d \tau = \int_0^\infty \langle
u_{k}(\cdot,T\tau), \Phi(\frac{k(T\tau)}{k(T)} \cdot, \tau)\rangle
d \tau \, .
\end{equation}
Assuming that \eqref{one} holds, the goal is to obtain a
nontrivial limit in \eqref{two} as $T \to \infty$. Then a natural
assumption is that
\begin{equation}
\label{three}
\lim_{T \to \infty} \frac{k(T \tau)}{k(T)} = p(\tau)
\end{equation}
exists for all $\tau \in \R^+$, since in this case
\[
\Phi \left( \frac{k(T\tau)}{k(T)} x, \tau \right) \, \to \Phi
(p(\tau)x, \tau )
\]
uniformly with all derivatives in $x$, boundedly in $\tau$, and
thus
\[
\int_0^\infty \langle  u_{T,k}(\cdot,\tau),
\Phi(\cdot,\tau)\rangle d \tau \to \int_0^\infty \langle U,\Phi
(p(\tau)\cdot, \tau ) \rangle\, d\tau.
\]

\medskip
A function $k$ that is eventually positive and for which
\eqref{three} holds for some function $p$, for all $\tau \in
(0,\infty)$ is called \emph{regularly varying} (\cite{bingham}).
It is known that in this case $p(\tau) = \tau^\alpha$ for some
$\alpha \in \R$ which is called the \emph{index}. An equivalent
condition is
\begin{equation}
\lim_{T \to \infty} \frac{k(T \tau)}{b(T)} = p_1(\tau)
\end{equation}
for some $b, \, p_1$ and all $\tau$ in a neighborhood of $\tau =
1$. In this case, necessarily $p_1(\tau) = C\tau^\alpha$ for some
$C = p_1(1) > 0$, and one may choose $b(t) = k(t)$.

Since $k$ is increasing b y assumption, $\alpha $ must be
non-negative. If \eqref{three} holds merely for $\tau$ in a set of
positive Lebesgue measure, it must already hold for all positive
$\tau$, and the limit is uniform on closed subintervals of
$(0,\infty)$ (see \cite{bingham}). Moreover, $\alpha$ can be
recovered from the limit (which always exists)
\[
\alpha = \lim_{t \to \infty} \frac {\log \, k(t)}{\log t}\, .
\]
Returning to \eqref{two} and assuming now that $k$ is regularly
varying with index $\alpha \ge 0$, one obtains
\begin{equation}
\label{four}
 \int_0^\infty \langle  u_{T,k}(\cdot,\tau),
 \Phi(\cdot,\tau)\rangle d \tau \, \to \int_0^\infty \langle U, \Phi_\tau(\cdot,\tau) \rangle
d \tau
\end{equation}
where $\Phi_\tau(x,\tau) = \Phi(x\tau^\alpha,\tau)$. If $U$ is a
locally integrable function, this means
\[
u_{T,k}(x,\tau) \sim \tau^{-n \alpha}U(x\tau^{-\alpha}) \, .
\]
in a suitable sense (e.g. pointwise a.e. in $(x,\tau)$), as $T \to
\infty$. All this proves the first statement of the following
proposition.

\begin{proposition}
Let $t \mapsto u(\cdot,t)$ be a measurable locally bounded family
of tempered distributions on $\R^n$. Let $k:[0,\infty) \to \R^+$
be continuous and increasing to $\infty$, and regularly varying
with index $\alpha \ge 0$. Define $u_k$ and $u_{T,k}$ as above.

\noindent
 a) If $U$ is a distribution on $\R^n$ such that
$u_k(\cdot,\tau) \, \to \, U$, then $u_{T,k}(\cdot,\tau) \, \to
\bar U(\cdot,\tau)$ for a.e. $\tau$ in $\cD'$ as $T \to \infty$,
where
\begin{equation}
\label{foura} \langle \bar U(\cdot,\tau), \varphi \rangle =
\langle U, \varphi( \tau^\alpha \cdot) \rangle\, .
\end{equation}
b) Assume that $U \in H^s(\R^n)$ for some $s \in \R$ and that
$u_{T,k}$ converges to $\bar U$, defined in\eqref{foura}, locally
uniformly in $\tau$ as an $H^s$ - valued function. Then also
\[ u_k(\cdot,t) \, \to \, U
\]
in $H^s(\R^n)$, as $t \to \infty$.

\noindent
c) Suppose $U \in L^r(\R^n)$ for some $r \in[1,\infty]$
and $u_k(\cdot,t) \, \to \, U$ in $L^r(\R^n)$. Set $w(x,t) =
k(t)^{-n} U(xk(t)^{-1})$, then
\[ \|u(\cdot,t) - w(\cdot,t)\|_{L^r} = o(k(t)^{n(r^{-1}-1)}) \, .
\]
Moreover, if $D^mU \in L^r(\R^n)$ and $D^m u_k(\cdot,t) \, \to \,
D^m U$ in $L^r(\R^n)$ for some partial derivative $D^m$ of order
$m$, then
\[ \|D^mu(\cdot,t) - D^mw(\cdot,t)\|_{L^r} = o(k(t)^{n(r^{-1}-1)-m}) \, .
\]
\end{proposition}

\medskip
\begin{proposition}
Let $k,\, l$ be regularly varying functions with index $\alpha \ge
0$ such that $\lim_{t \to \infty} \frac{l(t)}{k(t)} = C \in (0,
\infty)$, and let $U \in \cD'$. Suppose that $u_{T,k}(\cdot,t) \to
\bar U(\cdot,t)$, defined as in \eqref{foura}, locally uniformly
in some $H^s$ with $s \in \R$. Then $u_{T,l}(\cdot,t) \to \bar V
(\cdot,t)$ in $H^s$, locally uniformly in $t$, where for $\varphi
\in \cD$
\[
\langle \bar V(\cdot,\tau), \varphi \rangle = \langle U,
\varphi(C^{-1} \tau^\alpha \cdot) \rangle \, .
\]
Thus if $U$ is a function, then
\[
u_{T,l}(x,\tau) \sim C^n\tau^{-n \alpha}U(Cx\tau^{-\alpha})
\]
in a suitable sense, and the estimates of part c) of the previous
proposition hold.
\end{proposition}
The easy proofs of the remaining parts of Proposition 2.1 and of
Proposition 2.2 are left to the reader.

\medskip
If $u(\cdot,t) \to 0$ in some weak sense as $t \to \infty$, then
it is possible to find $k(\cdot)$, going to $\infty$, such that
$u_{T,k} \to 0$ as a distribution on $\R^n \times (0,\infty)$.
Similarly, if $u(\cdot,t) \in L^1$ for all $t$ and
$\int_{\R^n}u(\cdot,t) = C$ is constant,then one expects that
$u_{T,k}(\cdot,t) \to C\delta_0$ if $k(\cdot)$ goes to $\infty$
sufficiently rapidly. The limiting cases $U = 0$ and $U(\cdot,t) =
C\delta_0$ should be excluded and will be called {\it trivial}.
Non-trivial limiting distributions $U$ should be neither zero nor
supported on $\{0\}\times [0,\infty) \subset \R^n \times
[0,\infty)$.

\section{Identifying Asymptotic Limits}

In this section, I shall classify the types of limiting behavior
that are possible for distributional solutions of the partial
integrodifferential equation
\begin{equation}
\label{ide}
 u_t(x,t) = a_0 \Delta u(x,t) + \int_0^t a(t-s) \Delta
u(x,s) ds
\end{equation}
or more shortly $u_t = a_0 \Delta u + a*\Delta u$ for $x \in \R^n,
\, t >0$, with initial data $u(\cdot,0) = u_0$. It will be assumed
throughout that $a_0 \ge 0$, the integral kernel $a$ is bounded on
any set $[\epsilon,\infty)$ and integrable on $(0,1)$, and $u_0
\in L^1(\R^n)$. Let us begin by describing the limiting equations
and their solutions. For $\alpha
> 0$, the \emph{Mittag-Leffler} function $E_\alpha$ is defined as
\begin{equation}
E_\alpha(z) = \sum_{k=0}^\infty \frac{z^k}{\Gamma(1+\alpha k)} \,
,
\end{equation}
see \cite{bateman}. This is an entire function for any $\alpha >
0$. Special cases include $E_1(z) = e^z, \, E_2(z^2) = \cosh(z)$,
and $E_{1/2}(z) = e^{z^2} erfc(-z)$, where $erfc$ is the
complementary error function. For $0 < \alpha < 2$, $\alpha \neq
1$, there is an asymptotic expansion
\[
E_\alpha(z) = \sum_{n=1}^N\frac{z^{-n}}{\Gamma(1-\alpha n)} +
O(|z|^{-N-1})
\]
as $z \to \infty$ in a sector about the negative real axis, where
the reciprocal of the $\Gamma$ - function is extended as zero at
the poles of $\Gamma$. In particular, $E_\alpha$ is bounded on the
negative real axis for all $\alpha \le 2$; see \cite{bateman} for
details and other properties.

\medskip
For $\alpha > 0, \, \lambda \in \R$ the function $u(t) =
E_\alpha(-\lambda t^\alpha)$ is the solution of the scalar
integral equation
\begin{equation}
\label{scalar0}
 u(t) + \frac{\lambda}{\Gamma(\alpha)} \int_0^t
(t-s)^{\alpha-1}u(s) \, ds = 1
\end{equation}
as a direct calculation shows.

\medskip
Let $0 < \alpha \le 2$. Define $w_\alpha(\cdot, t)$ as the
tempered distribution on $\R^n$ whose spatial Fourier transform is
\begin{equation}
\label{defw} \hat w_\alpha(\xi,t) = E_\alpha(-|\xi|^2 t^\alpha)
\quad (\xi \in \R^n, \, t \ge 0) \,
\end{equation}
for $t > 0$.  If $\alpha < 2$, the asymptotic behavior of
$E_\alpha$ implies that $w_\alpha(\cdot,t) \in H^s(\R^n)$ for $s <
2-n/2$. In particular, if $\alpha < 2$, with the exception $\alpha
= 1$, then $w_\alpha(\cdot,t) \in L^2(\R^n)$ if and only if $n \le
3$ and $w_\alpha(\cdot,t) \in L^\infty(\R^n)$ if and only if $n =
1$. For $\alpha = 2$, one has $\hat w_2(\xi,t) = \cos(|\xi|t)$ and
thus $w_2(\cdot,t) \in H^s(\R^n)$ for all $s<-n/2$. If $\alpha =
1$, then $\hat w_1(\xi,t) = e^{-|\xi|^2 t}$, and $w_1$ is the
well-known fundamental solution of the heat equation. Also,
$w_\alpha(\cdot,t) \, \to \, \delta_0$ as $t \, \to \, 0$ for all
$\alpha$, in the sense of distributions.

\medskip
The distribution $w_\alpha$ solves the integrodifferential
equation
\begin{equation}
\label{inte0} w(\cdot,t) - \Delta
\left(\int_0^t\Gamma(\alpha)^{-1}(t-s)^{\alpha-1}w(\cdot,s) \, ds
\right) = \delta_0 \, .
\end{equation}
This follows immediately from \eqref{scalar0} and \eqref{defw}. If
$\alpha = 1$, this is the heat equation, and if $\alpha = 2$, this
is the wave equation. For $1 < \alpha < 2$, equation (12) can be
differentiated formally, resulting in the \emph{fractional heat
equation}
\begin{equation}
w_t(\cdot,t) = \Delta \left(\int_0^t
\Gamma(\alpha-1)^{-1}(t-s)^{\alpha-2}w(\cdot,s) \, ds \right)  .
\end{equation}
These equations have been studied in \cite{fujita, mainardi1,
hanyga1, mainardi0, mainardi2, schneiderwyss}.

\medskip

Returning to solutions of \eqref{ide}, let us assume that $t
\mapsto u(\cdot,t)$ is a function with values in the set of
tempered distributions such that for all test functions $\Phi \in
C_0^\infty(\R^n \times [0,\infty))$ the equation holds
\begin{eqnarray*}
\int_0^\infty \langle u(\cdot,s), \Phi_t(\cdot,s) + a_0 \Delta
\Phi(\cdot,s) + \int_t^\infty a(t-s) \Phi(\cdot,t) \, dt\rangle \,
ds
\\
+ \langle u(\cdot,0), \Phi(\cdot,0)\rangle  = 0 \, .
\end{eqnarray*}
With $A(t) = a_0 + \int_0^t a(s)$, one can write equivalently
\begin{equation}
\label{inte} u(\cdot,t) = u(\cdot,0) + \Delta \left( \int_0^t
A(t-s) u(\cdot,s) \, ds \right)
\end{equation}
in the sense of distributions, or with $ \Phi_1(x,s) = \int_s
^\infty A(t-s) \Phi(x,t) \, dt$
\begin{equation}
\label{testinte}
 \int_0^\infty \langle u(\cdot,s), \Phi(\cdot,s) - \Delta
\Phi_1(\cdot,s) \rangle \, ds = \langle u(\cdot,0), \int_0^\infty
\Phi(\cdot,t) dt \rangle .
\end{equation}

\medskip
Let us consider solutions of \eqref{ide} with the scaling
\[
u_{T,k}(x,\tau) = k(T)^n u(k(T),T \tau)
\]
introduced earlier, where $k$ is left unspecified for now. Set $K
= k(T)$, then $v = u_{T,k}$ is seen to be a distributional
solution of the problem
\begin{equation}
\label{scaledide}
 v = u_{0,K}+ \frac{T}{K^2}A_T*\Delta v
\end{equation}
where $u_{0,K} = K^n u_0(Kx)$ and $A_T(t) = A(Tt) = a_0 +
\int_0^{Tt} a(s) \, ds$. Also, the spatial Fourier transform $\hat
u_{T,k}$ solves
\begin{equation}
\label{ftinte}
\hat u_{T,k}(\xi,t) + |\xi|^2 \frac{T}{K^2}
A_T*\hat u_{T,k}(\xi,t) = \hat u_0(\frac{\xi}{K})
\end{equation}
in the sense of distributions. If one wishes to obtain a limiting
equation of a similar form, one is led to assume that there exist
functions $p, A_\infty$ such
\[
\lim_{T \to \infty} \frac{A_T(t)}{p(T)} = A_\infty(t) \, .
\]
Let us also assume that $A$ is eventually positive (not
necessarily bounded away from zeros). As explained in the previous
section, this implies that $A$ is regularly varying and that one
may choose $p(T) = cA(T) = cA_T(1), A_\infty(t) = c^{-1}t^\beta$
for some $\beta \in \R$ and any constant $c > 0$. If the kernel
$A_\infty$ is to be integrable at $t=0$, then one should require
$\beta > -1$. Since $A' = a$ was assumed to be bounded on
$(1,\infty)$, necessarily $\beta \le 1$. This motivates the main
assumption in the following result.

\medskip
\begin{theorem}
Let $u$ be a solution of \eqref{ide} in the sense described above,
and let $A$ be eventually positive and regularly varying with
index $\beta \in (-1,1]$. Assume that there exist a non-decreasing
function $k_0: [0,\infty) \, \to \, [0,\infty)$ with $k_0(\infty)
= \infty$, a sequence $T_n \to \infty$ and a non-trivial limiting
distribution $u_\infty$ on $\R^n \times [0,\infty)$,
$u_\infty(\cdot,t) \in H^s(\R^n)$ for a.e. $t$ for some fixed $s >
-\infty$ such that
\[
u_{T_n,k_0} \, \to \, u_\infty
\]
as $n \, \to \, \infty$, a.e. boundedly in $H^s(\R^n)$. Then one
may choose $k(t) = \sqrt{tA(t)\Gamma(1 + \beta)}$, and with this
choice
\[
u_{T,k} \, \to \, U_0 w_{1 + \beta} = u_\infty
\]
where $U_0 = \int_{\R^n} u_0(x) \, dx$ and $w_{1 + \beta}$ is the
distributional solution of the integrodifferential equation
\[
w(\cdot,t) - \Delta \left(\int_0^t\Gamma(1+\beta)^{-1}(t-s)^\beta
w(\cdot,s) \, ds \right) = \delta_0
\]
defined in \eqref{defw}.
\end{theorem}

\medskip
\begin{proof}
It should be noted that $A$ is regularly varying with index
$\beta$ iff $k$ is regularly varying with index
$\frac{1+\beta}{2}$. The assumptions for $A$ imply that
\[
\frac{A_T(t)}{A(T) \Gamma(1+\beta)} \, \to A_\infty(t) =
\frac{t^\beta}{\Gamma(1+\beta)} \]
 as $T \to \infty$. Let $\Phi
\in C^\infty_0(\R^n,\times (0,\infty))$ be a test function such
that
\begin{eqnarray*}
\langle\langle \Delta \Phi_\infty,u_\infty \rangle\rangle &\neq&
0\\
\langle\langle u_\infty, \Phi\rangle \rangle &\neq& U_0
\int_0^\infty  \Phi(0,t) \,dx dt
\end{eqnarray*}
where $\Phi_\infty(x,s) = \int_t^\infty A_\infty(t-s) \Phi(x,t) \,
dt$. This is possible since the limiting distribution is
non-trivial (not identically equal to zero, not supported on a
subset of $\R^n \times [0,\infty)$ and in some fixed $H^s$ for
a.e. $t$. Set $\Phi_n(x,s) = \int_t^\infty \frac
{A(T_n(t-s))}{A(T_n)\Gamma(1+\beta)} \Phi(x,t) \, dt$. Then
$\Phi_n \to \Phi_\infty$ together with all derivatives.

\medskip
Let us first show that $k_0$ may be replaced with $k$, i.e.
$u_{T_n,k} \, \to \, u_\infty $. Indeed, since
\[
\langle \langle\Phi, u_{T_n,k_0} \rangle\rangle= \int_0^\infty
\langle \Phi(\cdot,t), u_{0,k(T_n)}\rangle \, dt +
\frac{\Gamma(1+\beta) T_n A(T_n)}{k_0(T_n)^2} \langle\langle\Delta
\Phi_n, u_{T_n,k_0}\rangle\rangle
\]
and all terms except the fraction have limits as $n \to \infty$,
it follows that
\[
L = \lim_{n \to \infty}\frac{\Gamma(1+\beta) T_n
A(T_n)}{k_0(T_n)^2}
\]
exists. Of course, $L > 0$, and after replacing $k_0$ with
$\sqrt{L}k_0$, one may assume without loss of generality that $L =
1$. One can therefore replace $k_0$ with $k(t) =
\sqrt{\Gamma(1+\beta)tA(t)}$ and obtain that $u_{T_n,k} \to
u_\infty$.

\medskip
Observe next that \eqref{ftinte} now takes the form
\begin{equation}
\hat u_{T,k}(\xi,t) +  \frac{|\xi|^2}{A(T)\Gamma(1+\beta)}
A_T*\hat u_{T,k}(\xi,t)= \hat u_0(\frac{\xi}{K})
\end{equation}
for all $\xi$. Thus one can write
\[ \hat u_{T,k}(\xi,t) =
z_T(|\xi|^2,t)\hat u_0\left(\frac{\xi}{K}\right) \] where
$z_T(\rho,\cdot)$ solves the equation
\begin{equation}z_T(\rho,\cdot) + \frac{\rho}{A(T)\Gamma(1+\beta)} A_T*z_T(\rho,\cdot) = 1.
\end{equation}
Note that $z_T(\rho,t) = z(\lambda,Tt)$, where $z(\lambda,\cdot)$
solves
\begin{equation}
\label{scalar}
z(\lambda,\cdot) + \lambda A*z(\lambda,\cdot) = 1
\end{equation}
with $\lambda =\frac{\rho}{TA(T)\Gamma(1+\beta)}$.

By Lemma A.2, as $T \to \infty$, $z_T(\rho,t)$ converges to $E_{1+
\beta}(-\rho t^{1+\beta})$, locally uniformly in $\rho$ and $t \ge
0$, and thus $\hat u_{T,k}$ converges pointwise in $\xi$, locally
uniformly in $t$, to $\hat v(\xi,t) = U_0 E_{1+
\beta}(-|\xi|^2t^{1+\beta}) = U_0 w_{1+\beta}(\xi,t)$, that is, a
solution of the limiting equation
\[
\hat v(\xi,t)  + |\xi|^2 \int_0^t
\frac{(t-s)^\beta}{\Gamma(1+\beta)} \hat v(\xi,s) \,ds  = U_0\,  .
\]
But for the subsequence $T_n$, the limit is $\hat u_\infty$.
Therefore, convergence holds along the full sequence $T \to
\infty$, and the limit is $U_0\hat w_{1 + \beta}$. By Parseval's
identity, the theorem follows.
\end{proof}

\medskip
To sum up, the possibilities for limiting behavior identified in
this result are the following:

1. Behavior like the fundamental solution of the wave equation
($\beta = 1$), expected if e.g. $a(\cdot) \sim c >0$ or $a(t) \sim
\left( \log t\right)^m$ as $t \to \infty$ for some real number
$m$. In this case, $k(t) \sim \sqrt{c}t$ or $k(t) \sim t (\log
t)^{m/2}$.

2. Behavior like the fundamental solution of the heat equation
($\beta = 0$), expected if e.g. $a(\cdot)$ is integrable and
either $a \ge 0$ or $\int_0^\infty a(s) ds + a_0 > 0$, but also if
e.g. $a(t) \sim t^{-1}$. If $a$ is integrable, then $k(t) \sim
\sqrt{At}$, where $A = a_0 + \int_0^\infty a(s) ds$, while if e.g.
$a(t) \sim t^{-1}$, then $k(t) \sim \sqrt{t\log t}$.

3. Behavior like the fundamental solution of a fractional
integrodifferential equation of order $1 + \beta$. If  $0 < \beta
< 1$, this is expected if e.g. $a(t) \sim t^{\beta - 1}$. If $-1 <
\beta < 0$, this may occur if e.g. $a$ is negative and integrable,
$\int_0^\infty a(s) ds = -a_0$, and $a(t) \sim -t^{\beta - 1}$. In
each case, $k(t) \sim t^{(1+\beta)/2}$. Of course, the behavior of
$k$ may be modified by additional logarithmic factors also in this
case.

\section{Asymptotic Distributional Limits}

The purpose of this section is to demonstrate that under mild
assumptions, solutions of (7) do converge to limiting solutions
under the scaling $u \rightsquigarrow u_{T,k}$. I shall start by
stating a general existence result for strong solutions that is
essentially well-known for the \textsc{Hilbert} space case; see
\cite{pruss}.

\medskip
\begin{proposition}
Let $u_0 \in H^s(\R^n)$ for some $s \in \R$ and assume that
$$
\lim_{h \downarrow 0} (a_0 + \Re \tilde a(i \omega + h)) \ge 0
$$
for all $\omega \in \R$. Then there exists a unique function $u
\in C([0,\infty), H^s(\R^n))$ that solves \eqref{ide} in the sense
of distributions and for which $u(\cdot,0) = u_0$. The Fourier
transform $\hat u$ is given by
\begin{equation}
\label{rep}
 \hat u(\xi,t) = z(|\xi|^2,t) \hat u_0(\xi)
\end{equation} where $z(\lambda,\cdot)$ is the solution of
\eqref{scalar}.
\end{proposition}

\medskip
The condition for the kernel $a$ is equivalent to the requirement
that its cosine transform is bounded below by $-a_0$ in the sense
of measures. Equivalently, the measure $a_0 \delta_0 + a(t)dt$ is
required to be positive definite (see \cite{gripenberg}). We do
not require that $u_0 \in L^1$ and cannot assert that $u(\cdot,t)
\in L^1$ for $t>0$.

\medskip

\begin{proof}
Consider the integral equations \eqref{scalar} for $\lambda \ge
0$. By Lemma A.1, the estimate $|z(\lambda, t)| \le 1$ holds for
all $\lambda$ and $t$. Then define $u(\cdot,t)$ as in \eqref{rep}.
Using Parseval's identity and the bound for $z(\lambda,\cdot)$,
one obtains $\|u(\cdot,t)\|_{H^s} \le \|u_0\|_{H^s}$ for all $t$.
Also, $\hat u(\xi,t) \to \hat u_0(\xi)$ pointwise a.e., and by
construction
\[
\frac{\partial}{\partial t} \hat u(\xi,t) + |\xi|^2 \left( a_0
\hat u(\xi,t) + a*\hat u(\xi,t) \right) = 0
\]
for almost all $\xi$. Since $z(\lambda,\cdot)$ is continuous,
locally uniformly in $\lambda$, and uniformly bounded, $\hat
u(\xi,\cdot)$ is also continuous with values in $H^s(\R^n)$. Thus
$u$ is a distributional solution of \eqref{ide}. Uniqueness
follows by taking the Fourier transform of a solution in this
class and recognizing that it must have the form \eqref{rep}.
\end{proof}

\medskip
Such solutions converge to the limiting solutions identified in
the previous section under the scaling introduced there, if the
primitive $A$ of $a$ is regularly varying. As explained above,
these are natural conditions. From now on $u_0$ will always be
assumed to be integrable.

\medskip
\begin{theorem}
Assume that

1. the kernel satisfies $\lim_{h \downarrow 0} (a_0 + \Re \tilde
a(i \omega + h)) \ge 0$ for all $\omega \in \R$

2. the primitive $A(t) = a_0 + \int_0^t a(s) ds$ is eventually
positive and regularly varying with index $\beta \in (-1, 1]$.

Set $k(t) = \sqrt{t A(t) \Gamma(1+ \beta)}$, then for the solution
$u$ of \eqref{ide} found in the previous proposition
\begin{equation}
u_{T,k}(\cdot,t) \, \to \, U_0 u_\infty (\cdot,t)
\end{equation}
in $H^s(\R^n)$, if $s < -n/2$. Here $\hat u_\infty(\xi,t) = E_{1 +
\beta}(-|\xi|^2t^{1+ \beta})$ and $U_0 = \int_{\R^n} u_0(x) dx$.
The convergence is uniform on any interval $[c,d] \subset
(0,\infty)$.
\end{theorem}

\medskip
A few remarks can serve to put the result in perspective. First,
if $\beta = 1$, then $\hat u_\infty(\xi,t) = \cos(|\xi|t)$, and
thus $u_\infty(\cdot ,t) \notin H^{-n/2}(\R^n)$. Thus for a result
that covers the entire range $\beta \in (-1,1]$, one cannot expect
convergence in better spaces than $H^s(\R^n)$ with $s < -n/2$.
Also, conditions 1 and 2 in the above result are independent. For
example, the kernel $a(t) = \cos(t)$ with $A(t) = \sin(t)$ and
$\tilde  a = \frac{1}{2} \left(\delta_i + \delta_{-i}\right)$ in
the sense of measures on $i \R$ satisfies condition 1, but is not
regularly varying. The kernel $a(t) = 1 - e^{-t}$ has the
antiderivative $A(t) = a_0 + t-1+e^{-t}$ which is regularly
varying with $\beta = 1$, but since  $\Re \tilde a(i\omega) = -
\frac{1}{1+\omega^2}$ for $\omega \ne 0$, condition 1 is not
satisfied if $a_0 < 1$. Finally, it should be recalled that
$L^1(\R^n) \subset H^s(\R^n)$ for $s < -n/2$, but not for larger
$s$.

\medskip

\begin{proof}
Let $u$ be the distributional solution in
$C([0,\infty),H^s(\R^n))$ constructed in proposition 4.1. Then
$\hat u_{T,k}$ satisfies
\[
\hat u_{T,k}(\xi,t) + \frac{|\xi|^2}{A(T)\Gamma(1+\beta)} A_T*\hat
u_{T,k}(\xi,t)= \hat u_0 \left(\frac{\xi}{k(T)}\right)
\]
and therefore with $k(T) = \sqrt{TA(T)\Gamma(1+\beta)}$,
\[\hat
u_{T,k}(\xi,t) = z(\frac{|\xi|^2}{k^2(T)},Tt) \hat u_0
\left(\frac{\xi}{k(T)}\right)
\]
where $z(\lambda,\cdot)$ solves \eqref{scalar} and therefore
$v_\lambda(t) = z\left(
\frac{\lambda}{TA(T)\Gamma(1+\beta)},Tt\right)$ solves
\eqref{vscalar}. By Lemma A.2,  as $T \to \infty$,
\[
v_\lambda(t) = z(\frac{\lambda}{TA(T)\Gamma(1+\beta)},Tt) \to
E_{1+ \beta}(-\lambda t^{1+\beta})
\]
for all $\lambda > 0$. Consequently,
$$
\hat u_{T,k}(\xi,t) \to \hat u_0(0)
E_{1+\beta}(-|\xi|^2t^{1+\beta})
$$
pointwise for all $\xi$, locally uniformly in $t$.

\medskip
To prove that convergence holds in $H^s(\R^n)$ for $s < -n/2$, one
invokes again Lemma A.1 to deduce that $|z(\lambda,t)| \le 1$ and
therefore
$$
|\hat u_{T,k}(\xi,t)| \le |\hat u_0\left(\frac{\xi}{k(T)}\right)|
\le C
$$
for all $\xi$ and $T$. Then for $t > 0$
$$
\|u_{T,k}(\cdot,t) - U_0u_\infty(\cdot,t)\|^2_{H^s} = \int_{\R^n}
(1 +|\xi|^2)^s |\hat u_{T,k}(\xi,t) - U_0
E_{1+\beta}(-|\xi|^{2+2\beta}t)|^2 \, .
$$
Since $2s < -n$, Lebesgue's dominated convergence theorem implies
the conclusion.

\end{proof}

\section{Asymptotic Limits in $L^2$}
In this section, equation \eqref{ide} will be considered under the
scaling $u \rightsquigarrow u_{T,k}$, with the goal of proving
that limiting solutions are attained in
$L^\infty_{loc}(0,\infty;L^2(\R^n))$ or more generally
$L^\infty_{loc}(0,\infty;H^s(\R^n))$ with $s \ge 0$, provided the
initial data are in $L^2$ or $H^s$. The limiting solutions were
identified in section 2 and are unbounded in any $L^r(\R^n), \,
r>1$ as $t \to 0$. Thus one cannot expect uniform convergence up
to $t=0$. A convergence result in $L^2$ or in a better space
allows one to obtain the exact asymptotic behavior of solutions of
(7) in this space to leading order, by Proposition 2.1(c).
\smallskip

I only have a result for the case where $A$ is regularly varying
with index $\beta = 0$, i.e. $0 < a_0 + \int_0^\infty a(t) dt <
\infty$. The result is independent of the space dimension. In this
case the limiting equation is the heat equation. Note that for
$\beta \ne 0$, the limiting distributional solution is in
$L^2(\R^n)$ for $t
> 0$ if and only if $n \le 3$. Thus a result that holds for all
spatial dimensions cannot be expected if $\beta \ne 0$.

\medskip
\begin{theorem}
Assume that \newline 1. the initial data satisfy $u_0 \in
L^1(\R^n) \cap H^s(\R^n)$ for some $s \ge 0$
\newline
2. for some $\alpha > 0$, $\tilde a$ can be extended to the half
plane $\{s \in \C | \Re s \ge -\alpha\}$ and either $a_0 > 0$ and
$a_0 + \tilde a(i \omega)) > 0$, or  $\Re \tilde a(-\alpha + i
\omega) \ge 0$ for all $\omega$.
\newline
Set $k(t) = \sqrt{t A(t)}$, then for the solution $u$ of
\eqref{ide} found in Proposition 3.1
\begin{equation}
\label{lim}
 u_{T,k}(\cdot,t) \, \to \, U_0 u_\infty (\cdot,t)
\end{equation}
in $H^s(\R^n)$, where $U_0 = \int_{\R^n} u_0(x) dx$. Here $\hat
u_\infty(\xi,t) = e^{-|\xi|^2t}$, that is, $u_\infty$ is the
fundamental solution of the heat equation. The convergence is
uniform on any compact subinterval of $(0,\infty)$.
\end{theorem}
\begin{proof}
Note that assumption 2 implies that $0 < a_0 +\int_0^\infty a(t)dt
< \infty$.  Thus $A$ is regularly varying with index $\beta = 0$.
From the proof of Theorem 4.2 one sees that
\begin{equation}
\hat u_{T,k}(\xi,t) = z\left(\frac{|\xi|^2}{TA(T)},Tt\right) \hat
u_0 \left(\frac{\xi}{k(T)}\right) = z(\lambda,Tt) \hat u_0
\left(\frac{\xi}{k(T)}\right)
\end{equation}
with $\lambda = \frac{|\xi|^2}{TA(T)}$. By Lemma A.1, there are
therefore estimates, valid for all sufficiently large $T$,
\[
\left|z\left(\frac{|\xi|^2}{TA(T)},Tt\right)\right|\le 2
e^{-c|\xi|^2t}
\]
if $|\xi|^2 \le dT$ and
\[
\left|z\left(\frac{|\xi|^2}{TA(T)},Tt\right)\right|\le 2 e^{-cdTt}
\]
if $|\xi|^2 \ge dT$. Here $c = \frac{\epsilon}{2A_\infty},\, d =
\frac{A_\infty L}{2}$  are positive constants, and $\epsilon, \,
L$ are as in Lemma A.1. Now consider
\begin{eqnarray}
&&\|u_{T,k}(\cdot,t) - U_0u_\infty(\cdot,t)\|_{H^s}\nonumber\\
&=&\left(\int\left(1+|\xi|^{2s}\right)|z(\lambda ,Tt)\hat
u_0\left(\frac{\xi}{k(T)}\right) - U_0 e^{-|\xi|^2t}|^2 d \xi\right)^{1/2}\nonumber\\
&\le&\left(\int_{|\xi|\le R}\left(1+|\xi|^{2s}\right)|z(\lambda
,Tt)\hat u_0\left(\frac{\xi}{k(T)}\right) - U_0 e^{-|\xi|^2t}|^2 d
\xi\right)^{1/2} \nonumber\\
&\quad +& \left(\int_{|\xi|\ge
R}\left(1+|\xi|^{2s}\right)|z(\lambda ,Tt)\hat
u_0\left(\frac{\xi}{k(T)}\right)|^2 d \xi \right)^{1/2}\nonumber\\
 &\quad +&
\left(\int_{|\xi|\ge R}\left(1+|\xi|^{2s}\right)U_0^2
e^{-2|\xi|^2t} d \xi\right)^{1/2} \nonumber
\end{eqnarray}
for arbitrary $R$, where as before $\lambda =
\frac{|\xi|^2}{TA(T)}$. Then for given $\delta > 0$ and $[a,b]
\subset (0,\infty)$, choose $R$ so large that the last integral is
less than $\delta$, uniformly in $t \in [a,b]$, and then choose
$T$ large enough such that the first integral is also bounded by
$\delta$, uniformly in $t$. This is possible because
$z\left(\frac{|\xi|^2}{TA(T)},Tt\right) \to e^{-|\xi|^2t}$ as $T
\to \infty$ by Lemma A.2, locally uniformly in $\xi$, and because
$\hat u_0$ is continuous. Then if $dT \ge R^2$, the second
integral in the last expression can be estimated by
\[
\dots \le \left(\int_{|\xi|\ge
R}\left(1+|\xi|^{2s}\right)e^{-2cdTt}|\hat
u_0\left(\frac{\xi}{k(T)}\right)|^2 d \xi \right)^{1/2} \le
e^{-cdTt}k(T)^{n/2+s}\|u_0\|_{H^s} \] which is also smaller than
$\delta$, if $T$ is chosen sufficiently large, locally uniformly
in $t$. The proof of this theorem is therefore complete.
\end{proof}
Using Proposition 2.2, one obtains
\begin{corollary}
Under the assumptions of Theorem 5.1, the solution $u$ of
\eqref{ide} satisfies
\[
\|u(\cdot,t) - w(\cdot,t) \|_{H^s} = o\left(t^{-n/4}\right)
\]
where $w(x,t) = U_0\left(4A_\infty \pi t\right)^{-n/2} \exp
\left(-|x|^2/(4A_\infty t)\right)$ is the solution of the heat
equation
\[
w_t = A_\infty \Delta w, \quad w(\cdot,0) = U_0 \delta_0
\]
and $A_\infty = a_0 + \int_0^\infty a(s) \,d s$.
\end{corollary}

\section{Linear Viscoelasticity}

Consider a viscoelastic material with mass density $\rho = 1$
occupying all $\R^3$, and denote the displacement of a material
point at position $x$ and time $t$ by $u(x,t)$ and the velocity at
this point by $v(x,t) = \partial_tu(x,t)$. Let us assume that the
material is at rest for $t<0$ and prescribe an initial velocity
field $v(\cdot,0) = v_0$. The \textsc{Boltzmann} model for linear
isotropic homogeneous viscoelasticity (\cite{leitman}) leads to
the equations of motion
\begin{equation}
\label{visco} v_t = a_0\Delta v + \frac{a_0 + 2b_0}{3} \nabla
\nabla \cdot v + a*\Delta v + \frac{a+2b}{3}*\nabla \nabla \cdot v
\, .
\end{equation}
Here $a_0\ge 0, \, b_0 \in \R$, and $a, \, b$ are suitable
scalar-valued functions that describe the stress response of the
material under shear and compression, respectively.

\medskip
The reader should recall the well-known decomposition into
divergence free and gradient components, as follows. For $u \in
H^s(\R^3,\R^3)$, let $Pu, \, Qu \in H^s$ be defined by
\[
P\hat u(\xi) = \frac{\xi \xi^T}{|\xi|^2} \hat u(\xi), \quad Q\hat
u(\xi) = \left(1-\frac{\xi \xi^T}{|\xi|^2}\right) \hat u(\xi) \, .
\]
Then for an $H^s$-valued solution $v$ of \eqref{visco}, one
obtains that $p=Pv$ and $q=Qv$ satisfy the equations
\begin{eqnarray}
p_t &=& \beta_0 \Delta p + \beta*\Delta p \\
q_t &=& a_0 \Delta q + a*\Delta q
\end{eqnarray}
with $\beta_0 = \frac{4a_0+2b_0}{3}$ and $\beta(t) =
\frac{4a(t)+2b(t)}{3}$ and with initial data $p(\cdot,0) = p_0 =
Pv_0$ and $q(\cdot,0) = q_0 = Qv_0 = v_0-p_0$. In addition,
$\nabla \times p = 0$ and $\nabla \cdot q = 0$ in the sense of
distributions. Thus $p$ and $q$ satisfy scalar integrodifferential
equations, and as in Proposition 4.1 one obtains an existence
result together with a representation formula for the solution,
given next.

\medskip
\begin{proposition}
Let $v_0 \in H^s(\R^3,\R^3)$ for some $s \in \R$ and assume that
\begin{eqnarray}
\lim_{h \downarrow 0} (a_0 + \Re \tilde a(i \omega + h)) &\ge& 0
\\
\lim_{h \downarrow 0} (\beta_0 + \Re \tilde \beta(i \omega + h))
&\ge& 0
\end{eqnarray}
for all $\omega \in \R$. Then there exists a unique function $v =
p+q \in C([0,\infty), H^s(\R^3, \R^3))$ that solves \eqref{visco}
in the sense of distributions and for which $v(\cdot,0) = v_0$.
The Fourier transforms $\hat p, \, \hat q$ are given by
\begin{eqnarray}
\label{rep1}
\hat p(\xi,t) &=& z_1(|\xi|^2,t) \hat p_0(\xi) \\
\label{rep2}
 \hat q(\xi,t) &=& z(|\xi|^2,t) \hat q_0(\xi)
\end{eqnarray}
where $z(\lambda,\cdot)$ is the solution of \eqref{scalar} and
$z_1$ solves \eqref{scalar} with $a_0, \, a(\cdot)$ replaced by
$\beta_0, \, \beta(\cdot)$.
\end{proposition}

\medskip
Let us now consider the case where $a, \beta \in L^1$,
corresponding to a viscoelastic material with vanishing elastic
equilibrium response, i.e. a liquid. In this case, the limiting
behavior is expected to resemble the fundamental solution of the
compressible \textsc{Stokes} system
\begin{equation}
\label{stokes}
 w_t = A \Delta w + (B-A) \nabla \nabla \cdot w
\end{equation}
where $A = a_0 + \int_0^\infty a(s) \, ds$ and $B = \beta_0 +
\int_0^\infty \beta(s) \, ds$. The fundamental solution is known
to be the matrix valued function $W(x,t) = U(x,Bt) + V(x,At)$,
where
\[
\hat U(\xi,t) = \frac{\xi \xi^T}{|\xi|^2}e^{-|\xi|^2t}, \quad \hat
V(\xi,t) = \left(E - \frac{\xi \xi^T}{|\xi|^2}\right)e^{-|\xi|^2t}
\]
with $E$ denoting the identity matrix. In real terms, $U$ and $V$
can be expressed in terms of error functions (Kummer functions,
confluent hypergeometric functions), e.g.
\[
U_{ij}(x) = \partial_i \partial_j \left(\frac{1}{4\pi|x|}
Erf\left(\frac{|x|}{\sqrt{4t}}\right)\right) \, .
\]
A recent derivation of these fundamental solutions in a more
general situation may be found in \cite{thomann}. Note that these
functions are not integrable with respect to $x$, since their
Fourier transforms are not continuous at $\xi = 0$; indeed they
behave like $O(|x|^{-3}$ as $|x| \to \infty$ for fixed $t>0$, due
to well-known asymptotic results for the Kummer function. Using
the arguments that led to the proof of Theorem 5.1, one can now
describe the asymptotic behavior of solutions of \eqref{visco} in
terms of $U$ and $V$. For this purpose, let us assume the
following:

\medskip
\begin{itemize}
\item For some $s \ge 0$, $v_0 \in L^1(\R^3,\R^3) \cap H^s(\R^3,
\R^3)$.

\item For some $\alpha > 0$, the Laplace transforms $\tilde a, \,
\tilde \beta$ can be extended to the half plane $\{z \, | \, \Re z
\ge - \alpha \}$

\item Either $a_0 > 0$ and $a_0 + \tilde a(i \omega)) > 0$, or  $\Re
\tilde a(-\alpha + i \omega) \ge 0$ for all $\omega$.

\item The same assumption for $\beta_0$ and $\beta(\cdot)$.
\end{itemize}
\medskip
As in the previous section, one can then use the representation
formulae \eqref{rep1} and \eqref{rep2} together with the results
of Appendix A to prove the following result.

\begin{theorem}
Under these assumptions, the solution $v$ of \eqref{visco}
satisfies
\[
\|v(\cdot,t) - V_0^TU(\cdot,Bt) - V_0^TV(\cdot,At)\|_{H^s} =
o\left(t^{-n/4}\right)
\]
where $V_0 = \int_{\R^3}v_0(x) \, dx \in \R^3$, $U$ and $V$ are
the components of the fundamental solution of the compressible
\textsc{Stokes} system \eqref{stokes}, and
\[A = a_0 + \int_0^\infty a(s) ds , \quad B = \frac{4}{3}A +
\frac{ 2}{3}\left(b_0 + \int_0^\infty b(s) ds \right)\, .
\]
\end{theorem}
\medskip
The result shows that to leading order for large $t$, the solution
$v(\cdot,t)$ of the \textsc{Boltzmann} system \eqref{visco}
behaves like the solution of the \textsc{Stokes} system
\eqref{stokes} with distributional initial data $w(\cdot,0) = V_0
\delta_0$.

\appendix
\section{Scalar Integral Equations}

In the following, let us assume that $z:[0,\infty) \, \to \R$ is a
solution of the scalar integrodifferential equation
\begin{equation}
\label{scalara}
 z'(t) + \lambda \left( a_0 z(t) + a*z(t)\right) =
0, \, z(0) = 1
\end{equation}
where $a_0 \ge 0, \, \lambda \ge 0$, and $a \in
L^1_{loc}(0,\infty;\R), a \in L^\infty(1,\infty;\R)$. Let $\tilde
a$ be the Laplace transform of $a$, defined for $\Re{s} > 0$.
Recall that by Parseval's identity, $a_0 + \Re \tilde {a}(s) \ge
c$ for all $s$ in the right half plane, for some non-positive
constant $c$, if and only if
$$
\int_0^T u(t) \left(a_0 u(t) + a*u(t)\right) dt \ge c \int_0^T
|u(t)|^2dt
$$
for all real-valued square integrable functions $u$ and all $T$;
see \cite{gripenberg}.
\medskip

\begin{lemma}
1. If $\lim_{h \downarrow 0} (a_0 + \Re \tilde a(i \omega + h))
\ge 0$ for all $\omega \in \R$ then
\[|z(t)| \le 1
\]
for all $t$ and all $\lambda \ge 0$.
\newline
2. Assume that for some $\alpha > 0$, $\tilde a$ can be extended
to $\{s \in \C | \Re s \ge -\alpha\}$ and that either $a_0 > 0$
and  $\inf_{\omega \in \R} (a_0 + \tilde a(i \omega)) > 0$, or
that $\Re \tilde a(-\alpha + i \omega) \ge 0$ for all $\omega$.
Then there is a constant $\epsilon > 0$ such that for all $t
> 0$
\begin{equation}|z(t)| \le 2e^{-\epsilon \min(\lambda, 1) t} \, .
\end{equation}
\end{lemma}

\begin{proof}
Let us consider the more general equation
\begin{equation}
\label{scalarb}
 z'(t) + \lambda
\left(a_0 + a*z(t) \right) = f(t)
\end{equation}
where $f \in L^1_{loc}(0,\infty;\R)$. It will be shown that
\newline
a) under the assumptions of part 1, for $\lambda = 1$,
\begin{equation}
\label{est1}
 |z(t)| \le 1 + \int_0^t |f(s)|ds
\end{equation}
for all $t$,
\newline
b) under the first set of assumptions in part 2, with $f = 0$,
\begin{equation}
\label{est2} |z(t)| \le e^{-\delta t}
\end{equation} for all $t$,
\newline
c) under the second set of assumptions and with $f = 0$,
\begin{equation}
\label{est3} |z(t)| \le 2e^{-\delta t}
\end{equation}
for all $t$. Here $\delta = \lambda \min \{\epsilon, 1\}$, and
$\epsilon > 0$ depends only on $a_0$ and $a$. Together these
assertions imply the lemma.

\medskip
To prove part a), multiply \eqref{scalara} with $z(t)$ and
integrate over $[0,T]$, resulting in the identity
\[
\frac{1}{2}|z(T)|^2 + \int_0^T z(t)\left(a_0 z(t) +a*z(t) \right)
dt = \frac{1}{2} + \int_0^T z(t)f(t) dt \, .
\]
Since $a_0 + \tilde a(i \omega) \ge 0$, the integral on the left
is non-negative, and the inequality $\frac{1}{2}|z(T)|^2 \le
\frac{1}{2} + \int_0^T |z(t)||f(t)| dt$ follows for all $T$.
Bihari's theorem now implies \eqref{est1}. Evidently this estimate
is independent of $\lambda \ge 0$.

\medskip
For the proof of b), set $z_\delta(t) = e^{\delta t}z(t)$
and $a_\delta(t) = e^{\delta t}a(t)$, where $\delta > 0$ will be
fixed later. Then $z_\delta$ satisfies
\begin{equation}
\label{scalarc}
 z_\delta'(t) + \lambda\left((a_0
- \delta\lambda^{-1})z_\delta(t) + a_\delta * z_\delta(t)\right) =
0
\end{equation}
and  $z_\delta(0) = 1$.
 Note that $\tilde a_\delta(s) = \tilde
a(s-\delta)$ whenever $\delta < \alpha$. Since $\Re \tilde a(s)$
is bounded on any vertical line $\Re s = \beta$ with $\beta
> - \alpha$, harmonic to the right of any such line, and bounded
away from $0$ near $\Re s = -\alpha$ one can find $\epsilon > 0$
such that $a_0 + \Re \tilde a(-\epsilon + i \omega) \ge \epsilon$
for all $\omega \in \R$. Then also $a_0 + \Re \tilde a(s) \ge
\epsilon$ whenever $\Re s > -\epsilon$.

Let now $ \lambda \le  1$. Set $\delta = \lambda \epsilon$, then
obviously
\[
a_0 - \delta \lambda^{-1} + \Re \tilde a_\delta (i \omega) = a_0 -
\epsilon  + \Re \tilde a(-\delta +  i \omega) \ge 0 \, .
\]
If $\lambda > 1$, one sets $\delta = \epsilon$ and obtains
\[
a_0 - \delta \lambda^{-1} + \Re \tilde a_\delta (i \omega) \ge a_0
- \epsilon  + \Re \tilde a(-\epsilon +  i \omega) \ge 0 \, .
\]
Part a), applied to \eqref{scalarc}, implies the desired estimate
in both cases.

To prove c), note that
\begin{equation}
\label{ten}
 \Re \tilde a(i \omega) \ge \frac{c_1}{\omega^2 + \alpha^2}
\end{equation}
for all $\omega$, for some $c_1 > 0$, since $\Re \tilde a$ is
harmonic and positive for $\Re s \ge - \alpha$. Also, $\arg \tilde
a(z) = 0$ for $\Re z = 0$ and $\arg \tilde a(-\alpha + i \omega)
\ge -\frac{\pi}{2}$ for $\omega \ge 0$ by assumption. Therefore by
the maximum principle for harmonic functions, $\arg \tilde a(z)
\ge \arg(z+\alpha)$ for all $z$ with $\Im z \ge 0, \, \Re z \ge
-\alpha$. This implies that
\[
\Im \tilde a(i\omega) \ge -\alpha \omega  \Re \tilde a(i \omega)
\]
for all $\omega \ge 0$. Now let $b(t) = d e^{-\alpha t}$ with
\[ d = \min \{\frac{c_1}{2 \alpha^2}, \, \frac{1}{2\alpha}\}
\]
where $c_1$ is as in \eqref{ten}. Thus $b(0) = d, \, b' = - \alpha
b$, and $\tilde b(s) = \frac{d}{s+\alpha}$. Consider the kernel
\[
a_1(t) = a(t) + \lambda a*b(t) - \alpha b(t)
\]
with Laplace transform
\[\tilde a_1(s) = \tilde a(s) \left(1 +
\frac{\lambda d}{s+\alpha}\right) - \frac{\alpha d}{s+\alpha}
\]
for $0 \le \lambda \le 1$. Then $\tilde a_1$ is analytic for $\Re
s > - \alpha$, and for $s = i \omega, \omega \ge 0$ one has
\begin{eqnarray*}
\Re \tilde a_1(i \omega) &&= \Re \tilde a(i\omega)\left(1 +
\frac{\lambda d \alpha}{\omega^2 + \alpha^2} \right) + \Im
 \tilde a(i \omega) \frac{\lambda d \omega}{\omega^2 +
 \alpha^2} - \frac{\alpha^2 d}{\omega^2 + \alpha^2} \\
 &&> \frac{c_1}{2(\omega^2 + \alpha^2)} + \frac{1}{2}\Re \tilde a(i\omega) -
 \Re \tilde a(i \omega) \frac{\lambda d \alpha\omega^2}{\omega^2 +
 \alpha^2} - \frac{\alpha^2 d}{\omega^2 + \alpha^2} \\
 &&\ge 0
\end{eqnarray*}
uniformly in $\lambda \in [0,1]$, by the choice of $d$. In
addition, $\Re \tilde a_1(s)$ is bounded uniformly in $\lambda$
and $s$ on the strip $-\alpha/2 \le \Re s $. One can therefore
find a positive $\gamma < \alpha$ such that whenever $\Re s \ge -
\gamma$, then $\frac{d}{2} + \Re a_1(s) \ge 0$.

\medskip
Now the estimate \eqref{est3} can be proved for small $\lambda$,
say $\lambda \le \min\{1, \frac{2\gamma}{3d}\}$. Forming the
convolution of \eqref{scalara} with $\lambda b$ and adding the
result to \eqref{scalara}, one obtains the equation
\[
z'(t) + \lambda b*z'(t) + \lambda \left(a + \lambda b*a \right)*z
= 0
\]
or equivalently
\[
z'(t) + \lambda \left( d z(t) + a_1*z(t)\right) = \lambda b(t) =
\lambda d e^{-\alpha t}\, .
\]
where  $a_1$ is as above (depending also on $\lambda$). Set
$\delta = \frac{\lambda d}{2}$ and as before $z_\delta(t) =
e^{\delta t} z(t), \, a_{1,\delta}(t) = e^{\delta t}a_1(t)$. The
resulting equation for $z_\delta$ is
\[
z_\delta'(t) + \lambda\left(\frac{d}{2}z_\delta(t) +
 a_{1,\delta} * z_\delta(t)\right) = \lambda d e^{(\delta - \alpha) t}, \,
z_\delta(0) =1 \, .
\]
Since $\Re \tilde a_{1,\delta}(i \omega) = \Re \tilde a_1(-\delta
+ i \omega) \ge - \frac{d}{2}$, part a) of this proof implies that
\[
|z_\delta(t)| \le 1 + \int_0^\infty \lambda d e^{-\alpha t +
\lambda d t/2} = 1 + \frac{\lambda d}{\alpha - \lambda d/2} \le 2
\]
where the last inequality follows from $\lambda d \le
\frac{2}{3}\gamma$. This implies
\[
|z(t)| \le 2 e^{-\lambda \epsilon t}
\]
with $\epsilon = \frac{d}{2}$, whenever $\lambda \le \min
\{\frac{2\gamma}{3 d}, 1\}$. Reducing $\epsilon$ if necessary,
inequality \eqref{est3}is proved for $\lambda \le 1$.

\medskip
The proof for $\lambda \ge 1$ is similar. One considers the kernel
function
\[
a_2(t) = a(t) + a*b(t) - \lambda^{-1} \alpha b(t)
\]
with Laplace transform
\[
 \tilde a_2(s) = \tilde a(s) \left(1 + \frac{d}{s+\alpha}\right) -
\frac{\lambda^{-1} \alpha d}{s+\alpha} \, .
\]
Then by the same argument, $\frac{d}{2} + \Re a_2(s) \ge 0$
whenever $\Re s \ge - \kappa$, for some $\kappa
> 0$ that does not depend on $\lambda \ge 1$.
Then form the convolution of \eqref{scalara} with $b$ and add the
result to \eqref{scalara}. This yields the equation
\[
z'(t) + d z(t) + \lambda\left( a_2*z(t)\right) = b(t) = d
e^{-\alpha t}\, .
\]
with  $a_2(t) = a(t) + a*b(t) - \lambda^{-1}\alpha b(t)$. Set
$\delta = \min \{\kappa, \frac{d}{2}, \alpha - d\}$ and $z_\delta
= e^{\delta t}, \, a_{2,\delta} = e^{\delta t} a_2(t)$. Then
\[
z_\delta'(t) + (d - \delta)z_\delta(t) +\lambda\left(
 a_{2,\delta} * z_\delta(t)\right) = d e^{(\delta - \alpha) t}, \,
z_\delta(0) =1 \, .
\]
Since $\Re \tilde a_{2,\delta}(i \omega) = \Re \tilde a_2(-\delta
+ i \omega) \ge - \frac{d}{2}$, part a) again implies that
\[
|z_\delta(t)| \le 1 + \int_0^\infty d e^{(\delta -\alpha) t} = 1 +
\frac{d}{\alpha - \delta} \le 2\, ,
\]
where the last inequality follows from $\delta \le \alpha - d$.
Therefore,
\[
|z(t)| \le 2 e^{-\delta t} \le 2e^{-\epsilon t}
\]
if $\epsilon$ as chosen earlier or possibly lowered, whenever
$\lambda \ge 1$. The proof is now complete.
\end{proof}
\medskip
\begin{lemma}
Let $(A_n)_{n \ge 1}$ be a sequence in $L^1(0,T_0;\R)$ such that
$\|A_n - A_\infty\|_{L^1} \to 0$. For $\rho \ge 0$, let
$w_n(\rho,\cdot)$ be the solution of
\[w_n(\rho,t) + \rho A_n*w_n(\rho,t) = 1\]
for $0 \le t \le T_0$.  Then
\[w_n(\rho,\cdot) \to w_\infty(\rho,\cdot) \]
uniformly in $t \in [0,T_0]$, locally uniformly in $\rho$, where
\[
w_\infty(\rho,t) + \rho A_\infty*w_\infty(\rho,t) = 1 \, .
\]
In particular, asume that  $A$ is a regularly varying kernel with
index $\beta \in (-1,1]$ and eventually positive, and set $A_T(t)
= A(Tt)$. Then the solutions $v_\lambda(\cdot)$ of
\begin{equation}
\label{vscalar}
 v_\lambda(t) +
\frac{\lambda}{A(T)\Gamma(1+\beta)}A_T*v_\lambda(t) = 1
\end{equation}
converge uniformly in $t$ and locally uniformly in $\lambda$ to
$E_{1+\beta}(-\lambda t^{1+\beta})$.
\end{lemma}
\begin{proof}
The first assertion follows from a standard argument for Neumann
series. Since $\frac{A(Tt)}{A(T)\Gamma(1+\beta)} \to
\frac{t^\beta}{\Gamma(1+\beta)}$, pointwise in $t$ and also in
$L^1(0,T_0)$, the second assertion follows as well.
\end{proof}

\section{Two Examples}
Here are two explicit examples of integrodifferential equations of
the form \eqref{ide} for which the assumptions in the main results
are not satisfied and the conclusions fail as well.

\smallskip
First consider \eqref{ide} with the kernel $a(t) =
\cos(t)$ and $a_0 = 0$. Thus $\tilde a(s) = \frac{s}{s^2 +1}$, and
therefore this kernel is positive definite; $\Re \tilde a(s) \ge
0$ for all $s$ in the right half plane. Since $A(t) = \sin(t)$,
the kernel is not regularly varying. Taking the Laplace transform
with respect to $t$ and the Fourier transform with respect to $x$,
one obtains that the Fourier-Laplace transform solution $\hat
{\tilde u}$ satisfies
\[
s\hat {\tilde u} + \frac{|\xi^2|s}{s^2 +1} \hat {\tilde u} = \hat
u_0 \, .
\]
After solving for $\hat {\tilde u}$ and inverting the Laplace
transform, one obtains
\[
\hat u(\xi,t) = \left( \frac{1}{|\xi|^2+1} +
\frac{|\xi|^2}{|\xi|^2 + 1} \cos
\left(\sqrt{1+|\xi|^2}t\right)\right) \hat u_0(\xi) \, .
\]
Thus $u(x,t) = u_1(x) + u_2(x,t)$, where $u_1 -\Delta u_1 = u_0$
and $u_2$ solves the Klein-Gordon equation $u_{2,tt} + u_2 =
\Delta u_2$ with initial data $u_2(\cdot,0) = u_0 - u_1, \,
u_{2,t}(\cdot,0) = 0$. Locally in $x$, $u(\cdot,t) \to u_1$ as $t
\to \infty$, since the contributions from $u_2$ are radiated off
to infinity. There is a nontrivial time-asymptotic limit (attained
e.g. pointwise a.e. for sufficiently smooth initial data) that
depends on the initial data.

\smallskip
As a second example, consider \eqref{ide} with the kernel $a(t) =
-e^{-t}$ and $a_0 = 1$. Thus $a_0 + \tilde a(s) =  \frac{s}{s
+1}$, and this kernel is also positive definite. Here $A(t) =
e^{-t}$, and the kernel $A$ can be viewed as regularly varying
with index $\beta = -\infty$. As before, taking the Laplace
transform with respect to $t$ and the Fourier transform with
respect to $x$, one obtains that the Fourier-Laplace transform
solution $\hat {\tilde u}$ satisfies
\[
s\hat {\tilde u} + \frac{|\xi^2|s}{s +1} \hat {\tilde u} = \hat
u_0 \, .
\]
The equation can be solved for $\hat {\tilde u}$ and the Laplace
transform can be inverted, and the result is
\[
\hat u(\xi,t) = \left( \frac{1}{|\xi|^2+1} +
\frac{|\xi|^2}{|\xi|^2 + 1} e^{-(1+|\xi|^2)t}\right) \hat u_0(\xi)
\, .
\]
In this case therefore $u(x,t) = u_1(x) + u_2(x,t)$, where as
before $u_1 -\Delta u_1 = u_0$ and $u_2$ now solves the diffusion
equation  $u_{2,t} + u_2 = \Delta u_2$ with initial data
$u_2(\cdot,0) = u_0 - u_1$. Again, there a nontrivial
time-asymptoptic limit that depends on the initial data, namely
$u(\cdot,t) - u_1 = O(t^{-n/2}e^{-t})$ as $t \to \infty$,
uniformly in $x$.


\begin{thebibliography}{99}
\bibitem[1]{bingham} N. H. Bingham, C. M. Goldie, J. L. Teugels,
{\it Regular Variation}. Cambridge Univ. Press, Cambridge, 1989.

\bibitem[2]{barenblatt} G. I. Barenblatt, {\it Scaling, Self-Similarity, and Intermediate
Asymptotics}, Cambridge University Press, Cambridge, 1996.

\bibitem[3]{bateman} A. Erdelyi (ed.), {\it Higher Transcendental
Function, vol. 3}. McGraw-Hill, New York, 1955.

\bibitem[4]{dassios} G. Dassios, F. Zafiropoulos, Equipartition of
energy in linearized three-dimensional viscoelasticity. Quart.
Appl. Math. {\bf 48} (1990), 715 -- 730.

\bibitem[5]{eidelman}S. D. Eidelman, A. N. Kochubei, Cauchy problem
for fractional diffusion equation. J. Diff. Eq. {\bf 199} (2004),
211 -- 255.

\bibitem[6]{fujita} Y. Fujita, Integrodifferential
equation which interpolates the heat equation and the wave
equation.  Osaka J. Math.  {\bf 27} (1990),  309--321. II, ibid.,
797 -- 804.

\bibitem[7]{mainardi1} R. Gorenflo, Y. Luchko, F. Mainardi, Wright
functions as scale-invariant solutions of the diffusion-wave
equation. J. Comp. Appl. Math. {\bf 118} (2000),175 -- 191.

\bibitem[8]{gripenberg} G. Gripenberg, S. O. Londen, O. J. Staffans,
{\it Volterra Integral and Functional Equations}, Cambridge Univ.
Press, Cambridge, 1990.

\bibitem[9]{hanyga1} A. Hanyga, Multidimensional solutions of
time-fractional diffusion-wave equations. Proc. Royal Soc. London
A {\bf 458} (2002), 933 -- 957.

\bibitem[10]{hanyga2} A. Hanyga, E. Rok, Wave propagation in
micro-heterogeneous porous media: A model based on an
integro-differential wave equation. J. Acoustical Soc. America
{\bf 107} (2000), 2965 -- 2972.

\bibitem[11]{mainardi0} F. Mainardi, The fundamental solutions for
the fractional diffusion-wave equation. Appl. Math. Letters {\bf
9} (1996), 23 -- 28.

\bibitem[12]{mainardi2} F. Mainardi, R. Gorenflo, On Mittag-Leffler
type functions in fractional evolution processes. J. Comp. Appl.
Math. {\bf 118} (2000), 283 -- 299.

\bibitem[13]{leitman} M. J.  Leitman, G. C. Fisher, The linear
theory of viscoelasticity. In: S. Fl{\"u}gge (ed.), {\it Handbuch
der Physik VIa/3}, p. 1 -- 123, Springer, New York, Heidelberg,
Berlin 1972.

\bibitem[14]{pipkin} A.C. Pipkin, Asymptotic behavior of
viscoelastic waves. Quart. J. mech. Appl. Math. {\bf 41} (1988),
51 -- 64.

\bibitem[15]{pruss} J. Pr{\"{u}}ss,{\it Evolutionary Integral Equations and Applications},
Birkh{\"{a}}user Verlag, Basel, 1993.

\bibitem[16]{schneiderwyss} W. R. Schneider, W. Wyss, Fractional
Diffusion and Wave Equations. J. Math. Phys. {\bf 30} (1989), 134
-- 144.

\bibitem[17]{thomann} E. A. Thomann, R. B. Guenther, The fundamental solution
of the linearized Navier-Stokes equations for spinning bodies in
three spatial dimensions. J. Math. Fluid Mech. (2005), online.

\bibitem[18]{zaslavsky} G. M. Zaslavsky, Chaos, fractional
kinetics, and anomalous transport. Physics Reports {\bf 317}
(2000), 461 -- 580.

\end{thebibliography}
\end{document}